\documentclass[aop,MSNbibl,seceqn,dvips]{arximspdf}

\doi{10.1214/13-AOP875} 
\volume{42}
\issue{3}
\pubyear{2014}
\firstpage{994}
\lastpage{1019}

\makeatletter
\newtheorem{theorem}{Theorem}[section]
\newtheorem{cor}[theorem]{Corollary}
\newtheorem{lem}[theorem]{Lemma}
\newtheorem{prp}[theorem]{Proposition}
\newproclaim{defn}{Definition}[section]
\newproclaim{rem}{Remark}[section]
\newproclaim{exa}{Example}[section]
\def\B{\mathcal B}
\def\pp{\partial}
\def\d{\mathrm{d}}
\makeatother

\begin{document}
\begin{frontmatter}

\title{Integration by parts formula and shift Harnack inequality for stochastic equations}
\runtitle{Integration by parts and shift Harnack inequality}

\begin{aug}
\author[A]{\fnms{Feng-Yu} \snm{Wang}\corref{}\ead[label=e1]{wangfy@bnu.edu.cn}\ead[label=e2]{F.Y.Wang@swansea.ac.uk}\thanksref{t1}}
\runauthor{F.-Y. Wang} \affiliation{Beijing Normal University and
Swansea University}
\address[A]{School of Mathematical Sciences\\
Beijing Normal University\\
Beijing 100875\\
China\\
\printead{e1}\\
and\\
Department of Mathematics\\
Swansea University\\
Singleton Park, Swansea SA2 8PP\\
United Kingdom\\
\printead{e2}} 
\end{aug}
\thankstext{t1}{Supported in part by Lab. Math. Com. Sys., NNSFC (11131003),
SRFDP and the Fundamental Research Funds for the Central Universities.}

\received{\smonth{3} \syear{2012}} \revised{\smonth{2} \syear{2013}}

%
\begin{abstract}
A new coupling argument is introduced to establish Driver's integration
by parts formula and shift Harnack inequality. Unlike known coupling
methods where two marginal processes with different starting points are
constructed to move together as soon as possible, for the new-type
coupling the two marginal processes start from the same point but their
difference is aimed to reach a fixed quantity at a given time. Besides
the integration by parts formula, the new coupling method is also
efficient to imply the shift Harnack inequality. Differently from known
Harnack inequalities where the values of a reference function at
different points are compared, in the shift Harnack inequality the
reference function, rather than the initial point, is shifted. A~number
of applications of the integration by parts and shift Harnack
inequality are presented. The general results are illustrated by some
concrete models including the stochastic Hamiltonian system where the
associated diffusion process can be highly degenerate, delayed SDEs and
semi-linear SPDEs.
\end{abstract}

%

\end{frontmatter}

\section{Introduction}\label{sec1}

In stochastic analysis for diffusion processes, the Bismut formula
\cite{B} (also known as Bismut--Elworthy--Li formula due to \cite{EL})
and the integration by parts formula are two fundamental tools. Let,
for instance, $X(t)$ be the (nonexplosive) diffusion process generated
by an elliptic differential operator on a Riemannian manifold $M$, and
let $P_t$ be the associated Markov semigroup. For $x\in M$ and $U\in
T_x M$, the Bismut formula is of type
%
%
\begin{equation}
\label{1.1} \nabla_U P_t f (x)= \mathbf E \bigl\{f
\bigl(X^x(t) \bigr)M^x(t) \bigr\},\qquad f\in
\B_b(M), t>0,
\end{equation}
where $X^x(t)$ is the diffusion process starting at point $x$, $M^x(t)$
is a random variable independent of $f$ and $\nabla_U$ is the
directional derivative along $U$. When the curvature of the diffusion
operator is bounded below, this formula is available with $M^x(t)$
explicitly given by $U$ and the curvature operator. There exist a
number of applications of this formula, in particular, letting
$p_t(x,y)$ be the density (or heat kernel) of $P_t$ w.r.t. a nice
reference measure $\mu$, we have, formally,
\[
\nabla_U \log p_t(\cdot, y) (x)=\mathbf E
\bigl(M^x(t)\mid X^x(t)=y \bigr).
\]
From (\ref{1.1}) one may also derive gradient-entropy estimates of
$P_t$ and thus, the following Harnack inequality introduced in \cite
{W97} (see \cite{ATW09,GW}):
%
%
\begin{eqnarray}
\label{1.2} |P_t f|^p(x)\le P_t
|f|^p(y) \mathrm{e}^{C_p(t,x,y)},
\nonumber\\[-10pt]\\[-10pt]
\eqntext{t>0, p>1, x,y\in M,f\in B_b(M),}
\end{eqnarray}
where $C_p(t,x,y)$ is determined by moments of $M^\cdot(t)$ and thus,
independent of~$f$. This type of Harnack inequality is a powerful tool
in the study of contractivity properties, functional inequalities and
heat kernel estimates; see, for example, \cite{WY}~and references
within.

On the other hand, to characterize the derivative of $p_t(x,y)$ in $y$,
which is essentially different from that in $x$ when $P_t$ is not
symmetric w.r.t. $\mu$, we need to establish the following integration
by parts formula (see \cite{D}):
%
%
\begin{equation}
\label{1.3} \qquad P_t(\nabla_U f) (x)= \mathbf E \bigl\{ f
\bigl(X^x(t) \bigr) N^x(t) \bigr\},\qquad f\in
C_0^1(M), t>0, x\in M
\end{equation}
for a smooth vector field $U$ and some random variable $N^x(t)$.
Combining this formula with (\ref{1.1}), we are able to estimate the
commutator $\nabla P_t-P_t\nabla$ which is important in the study of
flow properties; see, for example, \cite{FLL}. Similar to
(\ref{1.1}), inequality (\ref{1.3}) can be used to derive a formula for
$\nabla_U \log p_t(x,\cdot)(y)$ and the shift Harnack inequality of
type
%
%
\begin{eqnarray}
\label{1.4} |P_t f|^p(x) \le P_t
\bigl(|f|^p\circ\exp [U] \bigr) (x)
\mathrm{e}^{C_p(t,x,y)},
\nonumber\\[-10pt]\\[-10pt]
\eqntext{t>0, p>1, x,y\in M,f\in B_b(M),}
\end{eqnarray}
where $\exp_x\dvtx  T_xM\rightarrow M, x\in M$, is the exponential map
on the Riemannian manifold. Differently from usual Harnack inequalities
like (\ref{1.2}), in (\ref{1.4}) the reference function $f$, rather
than the initial point, is shifted. This inequality will lead to
different heat kernel estimates from known ones implied by (\ref{1.2}).

Before moving on, let us make a brief comment concerning the study of
these two formulas. The Bismut formula (\ref{1.1}) has been widely
studied using both Malliavin calculus and coupling argument; cf.
\cite{WX,WZ,Z} and references within. Although (\ref{1.3}) also has
strong potential of applications, it is, however, much less known in
the literature due to the lack of efficient tools. To see that
(\ref{1.3}) is harder to derive than (\ref{1.1}), let us come back to
\cite{D} where an explicit version of (\ref{1.3}) is established for
the Brownian motion on a compact Riemannian manifold. Unlike the Bismut
formula which only relies on the Ricci curvature, Driver's integration
by parts formula involves both the Ricci curvature and its derivatives.
Therefore, one can imagine that in general (\ref{1.3}) is more
complicated (and hence harder to derive) than (\ref{1.1}).\vadjust{\goodbreak}

To establish the integration by parts formula and the corresponding
shift Harnack inequality in a general framework, in this paper we
propose a new coupling argument. In contrast to usual coupling
arguments where two marginal processes start from different points and
meet at some time (called the coupling time), for the new-type coupling
the marginal processes start from the same point, but their difference
reaches a fixed quantity at a given time.

In the next section, we will introduce some general results and
applications on the integration by parts formula and the shift Harnack
inequality using the new coupling method. The general result obtained
in Section~\ref{sec2} will be then applied in Section~\ref{sec3} to a
class of degenerate diffusion processes, in Section~\ref{sec4} to
delayed SDEs and in Section~\ref{sec5} to semi-linear SPDEs.

We remark that the model considered in Section~\ref{sec3} goes back to
the stochastic Hamiltonian system, for which the Bismut formula and the
Harnack inequalities have been investigated in \cite{GW,WZ,Z} by using
both coupling and Malliavin calculus. As will be shown in
Section~\ref{sec2.1} with a simple example of this model, for the study
of the integration by parts formula and the shift Harnack inequalities,
the Malliavin calculus can be less efficient than the new coupling
argument.

\section{Some general results}\label{sec2}
In Section~\ref{sec2.1} we first recall the argument of coupling by
change of measure introduced in \cite{ATW06,WX} for the Harnack
inequality and the Bismut formula, and then explain how can we modify
the coupling in order to derive the integration by parts formula and
the shift Harnack inequality, and introduce the Malliavin calculus for
the study of the integration by parts formula. In the second subsection
we present some applications of the integration by parts formula and
the shift Harnack inequalities to estimates of the heat kernel and its
derivatives.

For a measurable space $(E,\B)$, let $\B_b(E)$ be the class of all
bounded measurable functions on $E$, and $\B_b^+(E)$ the set of all
nonnegative elements in $\B_b(E)$. When $E$~is a topology space, we
always take $\B$ to be the Borel $\sigma$-field, and let $C_b(E)$
[resp., $C_0(E)$] be the set of all bounded (compactly supported)
continuous functions on $E$. If, moreover, $E$ is equipped with a
differential structure, for any $i\ge1$ let $C^i_b(E)$ be the set of
all elements in $C_b(E)$ with bounded continuous derivatives up to
order $i$, and let $C_0^i(E)= C_0(E)\cap C_b^i(E)$. Finally, a
contraction linear operator $P$ on $\B_b(E)$ is called a Markov
operator if it is positivity-preserving with $P1=1$.

\subsection{Integration by parts formula and shift Harnack inequality}\label{sec2.1}

\begin{defn}
Let $\mu$ and $\nu$ be two probability measures on a measurable space
$(E,\B)$, and let $X,Y$ be two $E$-valued random variables w.r.t. a
probability space $(\Omega,\mathcal F,\mathbf P)$.\vadjust{\goodbreak}
\begin{longlist}[(ii)]
\item[(i)] If the distribution of $X$ is $\mu$, while under another
probability measure $\mathbf Q$ on $(\Omega,\mathcal F)$ the
distribution of $Y$ is $\nu $, we call $(X,Y)$ a coupling by change
of measure for $\mu$ and $\nu$ with changed probability $\mathbf
Q$.

\item[(ii)] If $\mu$ and $\nu$ are distributions of two stochastic
    processes with path space $E$, a~coupling by change of measure for
    $\mu $ and $\nu$ is also called a coupling by change of measure for
    these processes. In this case $X$ and $Y$ are called the marginal
    processes of the coupling.
\end{longlist}
\end{defn}

Now, for fixed $T>0$, consider the path space $E^T:=E^{[0,T]}$ for some
$T>0$ equipped with the product $\sigma$-field $\B^T:=\B^{[0,T]}$. Let
$\{ P^x(A)\dvtx  x\in E, A\in\B^T\}$ be a transition probability such
that $P^x(\{\gamma\in E^T\dvtx  \gamma(0)=x\})=1, x\in E$. For any
$t\in [0,T]$, let $P_t(x,\cdot)=P^x(\{\gamma(t)\in\cdot\})$ be the
marginal distribution of $P^x$ at time~$t$. Then
\[
P_t f(x):=\int_E f(y)P_t(x,\d
y),\qquad f\in\B_b(E), x\in E
\]
gives rise to a family of Markov operators $(P_t)_{t\in[0,T]}$ on $\B
_b(E)$ with $P_0=I$.

In order to establish the Harnack inequality, for any two different
points \mbox{$x,y\in E$}, one constructs a coupling by change of
measure $(X,Y)$ for $P^x$ and $P^y$ with changed probability $\mathbf
Q=R\mathbf P$ such that $X(T)=Y(T)$. Then
\begin{eqnarray*}
\bigl|P_T f(y) \bigr|^p &=& \bigl|\mathbf E_\mathbf Qf
\bigl(Y(T) \bigr) \bigr|^p
\\
&=& \bigl|\mathbf E \bigl\{R f \bigl(X(T) \bigr) \bigr
\} \bigr|^p
\\
&\le& \bigl(\mathbf E|f|^p \bigl(X(T) \bigr) \bigr) \bigl(\mathbf
ER^{p/(p-1)} \bigr)^{p-1}
\\
&=& \bigl(P_T|f|^p(x)
\bigr) \bigl(\mathbf ER^{p/(p-1)} \bigr)^{p-1}.
\end{eqnarray*}
This implies a Harnack inequality of type (\ref{1.2}) if $\mathbf
ER^{p/(p-1)}<\infty$.

To establish the Bismut formula, let, for example, $E$ be a Banach
space, and $x,e\in E$. One constructs a family of couplings by change
of measure $(X^\varepsilon,X) $ for $P^{x+\varepsilon e}$ and $P^x$
with changed probability $\mathbf Q_\varepsilon:=R_\varepsilon\mathbf
P$ such that $X^\varepsilon(T)=X(T),  \varepsilon\in [0,1]$. Then, if
$N^x(T):=\frac{\d}{\d\varepsilon} R_\varepsilon |_{\varepsilon=0}$
exists in $L^1(\mathbf P)$, for any $f\in\B_b(E)$, we obtain
\begin{eqnarray*}
\nabla_eP_Tf(x)&=& \frac{\d}{\d\varepsilon} \mathbf E \bigl
\{R_\varepsilon f \bigl(X^\varepsilon(T) \bigr) \bigr\}\Big|_{\varepsilon=0}
\\
&=& \frac{\d}{\d\varepsilon} \mathbf E \bigl\{R_\varepsilon f \bigl(X(T) \bigr)
\bigr\}\Big|_{\varepsilon=0}
\\
&=&\mathbf E \bigl\{f \bigl(X(T) \bigr)N^x(T)\bigr\}.
\end{eqnarray*}
Therefore, the Bismut formula (\ref{1.1}) is derived.

On the other hand, for the integration by parts formula and shift
Harnack inequality we need to construct couplings with marginal
processes starting from the same point but their ``difference'' equals
to a fixed value at time $T$. For simplicity, below we only consider
$E$ being a Banach space. To extend the result to nonlinear spaces like
Riemannian manifolds, one would need to make proper modifications using
the geometric structure in place of the linear structure.

\begin{theorem}\label{T2.1}
Let $E$ be a Banach space and $x,e\in E$ and $T>0$ be fixed.
\begin{longlist}[(2)]
\item[(1)] For any coupling by change of measure $(X,Y)$ for $P^x$
and $P^x$ with changed probability $\mathbf Q= R\mathbf P$ such
that $ Y(T)=X(T)+e$, there holds the shift Harnack inequality
\[
\bigl|P_T f(x) \bigr|^p\le P_T \bigl
\{|f|^p(e+\cdot) \bigr\}(x) \bigl(\mathbf ER^{p/(p-1)}
\bigr)^{p-1},\qquad f\in\B_b(E)
\]
and the shift log-Harnack inequality
\[
P_T\log f(x)\le\log P_T \bigl\{f(e+\cdot) \bigr\}(x) +
\mathbf E(R\log R),\qquad f\in\B_b(E),f> 0.
\]

\item[(2)] Let $(X, X^\varepsilon), \varepsilon\in[0,1]$, be a
family of couplings by change of measure for $P^x$ and $P^x$
with changed probability $\mathbf Q
_\varepsilon=R_\varepsilon\mathbf P$ such that
\[
X^\varepsilon(T)= X(T)+\varepsilon e,\qquad \varepsilon\in(0,1].
\]
If $R_0=1$ and $N(T):=-\frac{\d}{\d\varepsilon} R_\varepsilon
|_{\varepsilon=0}$ exists in $L^1(\mathbf P )$, then
%
%
\begin{equation}
\label{IIT} P_T (\nabla_e f) (x) =\mathbf E \bigl\{ f
\bigl(X(T) \bigr)N(T) \bigr\},\qquad f,\nabla_e f\in
\B_b(E).
\end{equation}
\end{longlist}
\end{theorem}

\begin{pf}
The proof is similar to that introduced above for the Harnack
inequality and the Bismut formula.
\begin{longlist}[(2)]
\item[(1)] Note that $P_Tf(x)= \mathbf E\{Rf(Y(T))\}=\mathbf
E\{Rf(X(T)+e)\}$. We have
\begin{eqnarray*}
\bigl|P_T f(x) \bigr|^p &\le& \bigl(\mathbf E|f|^p \bigl(X(T)+e \bigr)
\bigr) \bigl(\mathbf ER^{p/(p-1)} \bigr)^{p-1}
\\
&=& P_T \bigl\{|f|^p(e+\cdot) \bigr\}(x) \bigl(\mathbf ER^{p/(p-1)}
\bigr)^{p-1}.
\end{eqnarray*}
Next, by the Young inequality (see \cite{ATW09}, Lemma~2.4), for
positive $f$ we have
\begin{eqnarray*}
P_T\log f(x)&=&\mathbf E \bigl\{R\log f \bigl(X(T)+e \bigr) \bigr\}
\\
&\le&\log\mathbf Ef \bigl(X(T)+e \bigr) +\mathbf E(R\log R)
\\
&=& \log
P_T \bigl\{f(e+\cdot ) \bigr\}(x) +\mathbf E(R\log R).
\end{eqnarray*}

\item[(2)] Noting that $P_T f(x)=
\mathbf E \{R_\varepsilon f(X^\varepsilon(T)) \} = \mathbf
E \{R_\varepsilon f(X(T) +\varepsilon e) \}$, we obtain
\[
0= \frac{\d}{\d\varepsilon} \mathbf E \bigl\{R_\varepsilon f \bigl(X(T) +\varepsilon e
\bigr) \bigr\}\Big|_{\varepsilon=0} = P_T(\nabla_e f) (x) -
\mathbf E \bigl\{f \bigl(X(T) \bigr) N(T) \bigr\},
\]
\end{longlist}
provided $R_0=1$ and $N(T):=-\frac{\d}{\d\varepsilon}
R_\varepsilon|_{\varepsilon=0}$ exists in $L^1(\mathbf P)$.\vadjust{\goodbreak}
\end{pf}

From Theorem \ref{T2.1} and its proof we see that the machinery of the
new coupling argument is very clear. So, in applications the key point
of the study lies in the construction of new type couplings.

Next, we explain how one can establish the integration by parts formula
using Malliavin calculus. Let, for example, $W:=(W(t))_{t\ge0}$ be the
cylindrical Brownian motion on an Hilbert space
$(H,\langle\cdot,\cdot\rangle,|\cdot|)$ w.r.t. a probability space
$(\Omega,\mathcal F,\mathbf P)$ with natural filtration $\{\mathcal
F_t\}_{t\ge0}$. Let
\[
H^1:= \biggl\{h\in C \bigl([0,T];H \bigr)\dvtx \|h
\|_{H^1}^2:=\int_0^T
\bigl|h'(s) \bigr|^2\,\d s<\infty \biggr\}
\]
be the Cameron--Martin space. For a measurable functional of $W$,
denoted by~$F(W)$, such that $\mathbf EF(W)^2<\infty$ and
\[
H^1\ni h \mapsto D_h F(W):= \lim_{\varepsilon\downarrow0}
\frac
{F(W+\varepsilon h)-F(W)}{\varepsilon}
\]
gives rise to a bounded linear operator. Then we write $F(W)\in
\mathcal D(D)$ and call $DF(W)$ the Malliavin gradient of $F(W)$. It is
well known that $(D,\mathcal D(D))$ is a~densely defined closed
operator on $L^2(\Omega,\mathcal F _T;\mathbf P)$; see, for example,
\cite{N}, Section~1.3. Let $(D^*,\mathcal D(D^*))$ be its adjoint
operator, which is also called the divergence operator.

\begin{theorem}\label{MA}
Let $H,W,D$ and $D^*$ be introduced above. Let $e\in H$ and
$X\in\mathcal D(D)$. If there exists $h\in\mathcal D(D^*)$ such that
$D_h X=e$, then
\[
\mathbf E(\nabla_e f) (X)= \mathbf E \bigl\{f(X)D^*h \bigr\},\qquad
f\in C_b^1(H).
\]
\end{theorem}

\begin{pf} Since $D_h X=e$, we have
\[
\mathbf E(\nabla_ef) (X)= \mathbf E(\nabla_{D_h X}f) (X)=
\mathbf E \bigl\{D_h f(X) \bigr\} = \mathbf E \bigl\{ f(X)D^*h \bigr\}.
\]\upqed
\end{pf}

Finally, as the integration by parts formula (\ref{IIT}) and by the
Young inequality (see \cite{ATW09}, Lemma~2.4) imply the
derivative-entropy inequality
\begin{eqnarray*}
\bigl|P_T(\nabla_ef) \bigr|&\le&\delta \bigl\{P_T(f
\log f)-(P_Tf)\log P_Tf \bigr\}
\\
&&{}+ \delta\log\mathbf E \biggl\{\exp \biggl[\frac
{|N(T)|}\delta \biggr]
\biggr\}P_Tf,\qquad \delta>0
\end{eqnarray*}
and the $L^2$-derivative inequality
\[
\bigl|P_T(\nabla_ef) \bigr|^2\le \bigl(\mathbf
EN(T)^2 \bigr) P_Tf^2,
\]
according to the following result it also implies shift Harnack
inequalities.

\begin{prp}\label{PIS}
Let $P$ be a Markov operator on $\B_b(E)$ for some Banach space $E$.
Let $e\in E$.
\begin{longlist}[(2)]
\item[(1)] Let $\delta_e\in(0,1)$ and $\beta_e\in C((\delta
_e,\infty)\times E; [0,\infty))$. Then
%
%
\begin{equation}
\label{W12ba} \qquad \bigl| P(\nabla_ef) \bigr|\le\delta \bigl\{ P(f\log f)-(Pf)\log
Pf \bigr\} +\beta_e(\delta,\cdot)Pf,\qquad \delta\ge
\delta_e
\end{equation}
holds for any positive $f\in C_b^1(E) $ if and only if
%
%
\begin{eqnarray}\label{W12bb}
(Pf)^p &\le& \bigl(P \bigl\{f^p(r e+\cdot)
\bigr\} \bigr)
\nonumber\\[-4pt]\\[-14pt]
&&{}\times\exp \biggl[\int_0^1
\frac{pr}{1+(p-1)s} \beta_e \biggl(\frac{p-1}{r+r
(p-1)s}, \cdot+sre \biggr)
\,\d s \biggr]\nonumber
\end{eqnarray}
holds for any positive $f\in\B_b(E), r\in(0,\frac1 {\delta_e})$ and
$p\ge \frac1{1-r\delta_e}$.

\item[(2)] Let $C\ge0$ be a constant. Then
%
%
\begin{equation}
\label{L2G} \bigl|P(\nabla_e f) \bigr|^2\le CP f^2,
\qquad f\in C_b^1(E),f\ge0
\end{equation}
is equivalent to
%
%
\begin{equation}
\label{SHT} Pf \le P \bigl\{f(\alpha e+\cdot) \bigr\} +|\alpha |
\sqrt{CPf^2},\qquad \alpha\in \mathcal R, f \in\B_b^+(E).
\end{equation}
\end{longlist}
\end{prp}

\begin{pf}
The proof of (1) is similar to that of \cite{GW}, Proposition 4.1,
while (2)~is comparable to \cite{W12b}, Proposition~1.3.
\begin{longlist}[(2)]
\item[(1)] Let $\beta(s)= 1+(p-1)s, s\in[0,1]$. By the monotone class
    theorem, it suffices to prove for $f\in C_b^1(E)$. Since
    $\frac{p-1}{r\beta (s)}\ge \delta_e$ for $p\ge\frac1
    {1-r\delta_e}$, it follows from (\ref{W12ba}) that
\begin{eqnarray*}
&&\frac{\d}{\d s} \log \bigl(P\bigl\{f^{\beta(s)}(sre +\cdot)\bigr\}(x)
\bigr)^{p/\beta(s)}
\\
&&\qquad = \frac1 {\beta(s)^2P\bigl\{f^{\beta(s)}(sre+\cdot)
\bigr\}(x)}
\\
&&\quad\qquad\times{} \bigl( p(p-1) \bigl[P\bigl\{\bigl(f^{\beta(s)}\log
f^{\beta(s)}\bigr) (sre+\cdot)\bigr\}
\\
&&\hspace*{88pt}{}-\bigl(P\bigl\{f^{\beta(s)}(sre+\cdot)\bigr\}\bigr)\log P
\bigl\{f^{\beta (s)}(sre+\cdot )\bigr\} \bigr]
\\
&&\hspace*{173pt}{}+ prP\bigl\{\nabla_ef^{\beta(s)}(sre+
\cdot)\bigr\} \bigr) (x)
\\
&&\qquad{}\ge- \frac{rp}{\beta(s)}\beta_e \biggl(\frac{p-1}{r\beta
(s)}, x+sre
\biggr),\qquad s\in[0,1].
\end{eqnarray*}
Taking the integral over $[0,1]$ w.r.t. $\d s$ we prove (\ref{W12bb}).

Next, let $z,e\in E$ be fixed, and assume that $P(\nabla_e f)(z)\ge0$
(otherwise, simply use $-e$ to replace $e$). Then (\ref{W12bb}) with
$p= 1+\delta_e r$ implies that
\begin{eqnarray*}
&&\delta\bigl\{(Pf)\log Pf\bigr\} (z) +\bigl| P(\nabla_e
f)\bigr|(z)
\\
&&\qquad = \limsup_{r\rightarrow0}
\frac{(P\{f(re+\cdot)\})^{1+\delta r}(z)-Pf(z)}{r}
\\
&&\qquad\le\limsup_{r\rightarrow0} \frac1 r
\biggl\{\bigl(Pf^{1+\delta r}\bigr)(z)
\\
&&\hspace*{77pt}{}\times \exp\biggl[ \int _0^{1}
\frac{(1+\delta r)r}{1+\delta rs}\beta_e\biggl(\frac {\delta}{1+\delta
r s}, \gamma (r)\biggr)\,\d r\biggr]-Pf(z)\biggr\}
\\
&&\qquad= \delta P(f\log f)(z) +\beta_e(\delta) Pf(z).
\end{eqnarray*}%
Therefore, (\ref{W12ba}) holds.

\item[(2)] Let $r>0$. For nonnegative $f\in C_b^1(E)$, (\ref{L2G})
    implies that
\begin{eqnarray*}
&& \frac{\d}{\d s} P\biggl\{ \frac{f}{1+srf}\bigl(\alpha(1-s)e+\cdot\bigr)\biggr\}
\\
&&\qquad = -P\biggl\{\frac{rf^2}{1+srf}\bigl(\alpha(1-s)e+\cdot\bigr)\biggr\}-\alpha
P\biggl\{\nabla_e\biggl(\frac f{1+srf}\biggr)\bigl(\alpha(1-s)e+\cdot\bigr)\biggr\}
\\
&&\qquad \le-r P\biggl\{\frac{f^2}{1+srf}\bigl(\alpha(1-s)e+\cdot\bigr)\biggr\}
\\
&&\quad\qquad{}
+|\alpha|\biggl(CP\biggl\{\frac{f^2}{(1+srf)^2}\bigl(\alpha(1-s)e+\cdot\bigr)\biggr\}\biggr)^{1/2}
\\
&&\qquad \le\frac{\alpha^2 C}{4 r}.
\end{eqnarray*}
Noting that
\[
\frac{f}{1+rf}=f-\frac{rf^2}{1+rf}\ge f-rf^2,
\]
we obtain
\[
Pf\le P\bigl\{f(\alpha e+\cdot)\bigr\} + r Pf^2 +\frac{\alpha^2 C}{4 r},\qquad  r>0.
\]
Minimizing the right-hand side in $r>0$, we prove (\ref{SHT}).

On the other hand, let $x\in E$. Without loss of generality we assume
that $P(\nabla_e f)(x)\le0$, otherwise it suffices to replace $e$ by
$-e$. Then (\ref{SHT}) implies that
\[
\bigl|P(\nabla_e f)(x)\bigr|=\lim_{\alpha\downarrow0} \frac{Pf(x)- P\{
f(\alpha e+\cdot)\}
(x)}{\alpha} \le\sqrt{CPf^2(x)}.
\]
\end{longlist}
Therefore, (\ref{L2G}) holds.
\end{pf}

To conclude this section, we would like to compare the new coupling
argument with known coupling arguments and the Malliavin calculus, from
which we see that the study of the integration by parts formula and the
shift Harnack inequality is, in general, more difficult than that of
the Bismut formula and the Harnack inequality.

First, when a strong Markov process is concerned, for a usual coupling
$(X(t), Y(t))$ one may ask that the two marginal processes move
together after the coupling time, so that to ensure $X(T)=Y(T)$, one
only has to confirm that the coupling time is not larger than the given
time $T$. But for the new coupling argument, we have to prove that at
time $T$, the difference of the marginal processes equals to a fixed
quantity, which cannot be ensured, even if the difference already
reached this quantity at a (random) time before $T$. From this we see
that construction of a new-type coupling is, in general, more difficult
than that of a~usual coupling.

Second, it is well known that the Malliavin calculus is a very
efficient tool to establish Bismut-type formulas. To see the difficulty
for deriving the integration by parts formula using Malliavin calculus,
we look at a simple example of the model considered in
Section~\ref{sec3}, that is, $(X(t),Y(t))$ is the solution to the
following degenerate stochastic equation on $\mathcal R^2$:
%
%
\begin{equation}
\label{EE} \cases{\d X(t)= Y(t)\,\d t, \vspace*{2pt}
\cr
\d Y(t)= \d W(t) +Z
\bigl(X(t),Y(t) \bigr)\,\d t,}
\end{equation}
where $W(t)$ is the one-dimensional Brownian motion and $Z\in
C^1_b(\mathcal R ^2)$. For this model the Bismut formula and Harnack
inequalities can be easily derived from both the coupling method and
Malliavin calculus; see \cite{GW,WZ,Z}. We now explain how can one
establish the integration by parts formula using Malliavin calculus.
For fixed $T>0$ and, for example, $e=(0,1)$, to derive the integration
by parts formula for the derivative along $e$ using Theorem \ref{MA},
one needs to find $h\in\mathcal D(D^*)$ such that
%
%
\begin{equation}
\label{DH} D_h \bigl(X(T),Y(T) \bigr)=e.
\end{equation}
To search for such an element $h$, we note that (\ref{EE}) implies
\[
\d \bigl(D_h X(t), D_h Y(t) \bigr)=
\bigl(0,h'(t) \bigr)\,\d t +G(t) \pmatrix{D_hX(t)
\cr
D_h Y(t)}\,\d t
\]
and
\[
D_hX(0)=D_h Y(0)=0,
\]
where
\[
G(t):= \pmatrix{0 &1
\cr
Z' \bigl(\cdot,Y(t) \bigr) \bigl(X(t)
\bigr) &Z' \bigl(X(t),\cdot \bigr) \bigl(Y(t) \bigr)}.
\]
Then, (\ref{DH}) is equivalent to
\[
\int_0^T\mathrm{e}^{\int_t^TG(s)\,\d s} \pmatrix{0
\cr
h'(t)}\,\d t = (0,1).
\]
It is, however, very hard to solve $h$ from this equation for general
$Z\in C_b^1(\mathcal R^2)$. On the other hand, we will see in
Section~\ref{sec3} that the coupling argument we proposed above is much
more convenient for deriving the integration parts formula for this
example.

\subsection{Applications}\label{sec2.2}

We first consider $E=\mathcal R^d$ for some $d\ge1$, and to estimate
the density w.r.t. the Lebesgue measure for distributions and Markov
operators using integration by parts formulas and shift Harnack
inequalities.

\begin{theorem}\label{T4.1}
Let $X$ be a random variable on $\mathcal R^d$ such that for some $N\in
L^2(\Omega\rightarrow\mathcal R^d;\mathbf P)$
%
%
\begin{equation}
\label{4.1} \mathbf E(\nabla f) (X)= \mathbf E \bigl\{ f(X)N \bigr\},\qquad f\in
C^1_b \bigl(\mathcal R^d \bigr).
\end{equation}
\begin{longlist}[(2)]
\item[(1)] The distribution $\mathbf P_X$ of $X$ has a density
$\rho$ w.r.t. the Lebesgue measure, which satisfies
%
%
\begin{equation}
\label{DD0}\nabla\log\rho(x)= -\mathbf E(N\mid X=x),\qquad \mathbf
P_X\mbox{-a.s.}
\end{equation}
Consequently, for any $e\in\mathcal R^d$ and any convex positive
function $H$,
\[
\int_{\mathcal R^d} \bigl\{H \bigl( |\nabla_e\log\rho | \bigr)\rho
\bigr\} (x)\,\d x \le\mathbf EH \bigl( \bigl|\langle e,N\rangle \bigr| \bigr).
\]

\item[(2)] For any $U\in C_0^1(\mathcal R^d;\mathcal R^d)$,
\[
\mathbf E(\nabla_Uf) (X)= \mathbf E \bigl\{f(X) \bigl( \bigl\langle
U(X),N \bigr\rangle -(\operatorname{div} U) (X) \bigr) \bigr\},\qquad f\in
C^1 \bigl(\mathcal R^d \bigr).
\]
\end{longlist}
\end{theorem}

\begin{pf}
(1) We first observe that if $\mathbf P_X$ has density $\rho$, then for
any $f\in C_0^1(\mathcal R^d)$,
\begin{eqnarray*}
\int_{\mathcal R^d} \bigl\{\rho(x)\nabla f(x) \bigr\}\,\d x &=&
\mathbf E (\nabla f) (X)
\\
&=& \mathbf E \bigl\{f(X)\mathbf E(N\mid X) \bigr\}
\\
&=&\int_{\mathcal R^d}
\bigl\{ f(x)\mathbf E(N\mid X=x) \bigr\}\mathbf P_X(\d x).
\end{eqnarray*}
This implies (\ref{DD0}). To prove the existence of $\rho$, let
$\rho_n$ be the distribution density function of $X_n:=X+\frac\zeta n,
n\ge1$, where $\zeta$ is the standard Gaussian random variable on
$\mathcal R^d$ independent of $X$ and $N$. It follows from (\ref{4.1})
that
\[
\mathbf E(\nabla f) (X_n) =\mathbf E \bigl\{\nabla f(\zeta/n+\cdot)
\bigr\}(X) = \mathbf E \bigl\{f(X_n)N \bigr\}.
\]
Then
\[
4 \int_{\mathcal R^d}|\nabla\sqrt{\rho_n}|^2(x)
\,\d x= \mathbf E|\nabla\rho_n|^2(X_n) \le
\mathbf EN^2 <\infty.
\]
So, the sequence $\{\sqrt{\rho_n}\}_{n\ge1}$ is bounded in
$W^{2,1}(\mathcal R ^d;\d x)$. Thus, up to a subsequence,
$\sqrt{\rho_n}\rightarrow \sqrt{\rho}$ in $L^2_{\mathrm{loc}}(\d x)$
for some nonnegative function $\rho$. On the other hand, we have
$\rho_n(x)\,\d x \rightarrow\mathbf P_X(\d x)$ weakly. Therefore,
$\mathbf P_X(\d x)=\rho(x)\,\d x$.

(2) As for the second assertion, noting that for $U=\sum_{i=1}^d
U_i\pp_i $ one has
\[
\nabla_Uf =\sum_{i=1}^d
\pp_i (U_if) - f\operatorname{div} U,
\]
it follows from (\ref{4.1}) that
\begin{eqnarray*}
\mathbf E(\nabla_Uf) (X)&=& \sum_{i=1}^d
\mathbf E \bigl\{\pp _i(U_if) (X) \bigr\}- \mathbf E\{f
\operatorname{div} U\}(X)
\\
&=&\sum_{i=1}^d \mathbf E \bigl
\{(U_if) (X)N_i \bigr\}- \mathbf E\{f\operatorname{div} U
\}
\\
&=&\mathbf E \bigl\{f(X) \bigl( \bigl\langle U(X),N \bigr\rangle-(
\operatorname{div} U) (X) \bigr) \bigr\}.
\end{eqnarray*}\upqed
\end{pf}

Next, we consider applications of a general version of the shift
Harnack. Let $P(x,\d y)$ be a transition probability on a Banach space
$E$. Let
\[
Pf(x)=\int_{\mathcal R^{d}} f(y)P(x,\d y),\qquad f\in\B_b
\bigl(\mathcal R^{d} \bigr)
\]
be the associated Markov operator. Let $\Phi\dvtx [0,\infty)\rightarrow
[0,\infty)$ be a strictly increasing and convex continuous function.
Consider the shift Harnack inequality
%
%
\begin{equation}
\label{PH} \Phi \bigl(Pf(x) \bigr)\le P \bigl\{\Phi\circ f(e+\cdot ) \bigr\}(x)
\mathrm{e}^{C_\Phi(x,e)},\qquad f\in\B_b^+(E)
\end{equation}
for some $x,e \in E$ and constant $C_\Phi(x,e)\ge0$. Obviously, if
$\Phi(r)=r^p$ for some $p>1$, then this inequality reduces to the shift
Harnack inequality with power $p$, while when $\Phi(r)=\mathrm{e}^r$,
it becomes the log shift Harnack inequality.

\begin{theorem}\label{T4.2} Let $P$ be given above and satisfy (\ref{PH})
for all $x,e\in E:=\mathcal R^d$ and some nonnegative measurable
function $C_\Phi$ on $\mathcal R^d\times\mathcal R^d$. Then
%
%
\begin{equation}
\label{AA0} \qquad\sup_{f\in\B_b^+(\mathcal R^d), \int
_{\mathcal R^d}
\Phi\circ f(x)\,\d x\le1} \Phi(Pf) (x)\le\frac1 {\int
_{\mathcal R^d} \mathrm{e}^{-C_\Phi(x,e)}\,\d e},\qquad x\in\mathcal
R^d.
\end{equation}
Consequently:
\begin{longlist}[(2)]
\item[(1)] If $\Phi(0)=0$, then $P$ has a transition density
$\varrho (x,y)$ w.r.t. the Lebesgue measure such that
%
%
\begin{equation}
\label{AA1} \int_{\mathcal R^d} \varrho(x,y)\Phi ^{-1}
\bigl( \varrho(x,y) \bigr)\,\d y\le\Phi ^{-1} \biggl(\frac1 {\int
_{\mathcal R^d} \mathrm{e}^{-C_\Phi
(x,e)}\,\d e} \biggr).
\end{equation}

\item[(2)] If $\Phi(r)=r^p$ for some $p>1$, then
%
%
\begin{equation}
\label{AA} \int_{\mathcal R^d}\varrho(x,y)^{p/(p-1)}\,\d y\le
\frac1 { (\int_{\mathcal R
^d} \mathrm{e}^{-C_\Phi(x,e)}\,\d e
)^{1/(p-1)}}.
\end{equation}
\end{longlist}
\end{theorem}

\begin{pf}
Let $f\in\B_b^+(\mathcal R^d)$ such that $\int_{\mathcal R^d}
\Phi(f)(x)\,\d x\le1$. By (\ref{PH}) we have
\[
\Phi(Pf) (x) \mathrm{e}^{-C_\Phi(x,e)}\le P \bigl\{\Phi\circ f (e+\cdot) \bigr\}
(x)=\int_{\mathcal R
^d} \Phi\circ f(y+e)P(x, \d y).
\]
Integrating both sides w.r.t. $\d e$ and noting that $\int_{\mathcal
R^d}\Phi \circ f(y+e)\,\d e =\int_{\mathcal R^d} \Phi\circ f(e)\,\d
e\le1$, we obtain
\[
\Phi(Pf) (x) \int_{\mathcal R^d} \mathrm{e}^{-C_\Phi(x,e)}\,\d e\le1.
\]
This implies (\ref{AA0}). When $\Phi(0)=0$, (\ref{AA0}) implies that
%
%
\begin{equation}
\label{AA2}\sup_{f\in\B_b^+(\mathcal R^d), \int
_{\mathcal R^d} \Phi\circ f(x)\,\d
x\le1} Pf (x)\le\Phi^{-1} \biggl(
\frac1 {\int_{\mathcal R^d} \mathrm{e}^{-C_\Phi
(x,e)}\,\d e} \biggr)<
\infty
\end{equation}
since by the strictly increasing and convex properties we have $\Phi
(r)\uparrow\infty$ as $r\uparrow\infty$. Now, for any Lebesgue-null set
$A$, taking $f_n=n1_A$ we obtain from $\Phi(0)=0$ that
\[
\int_{\mathcal R^d} \Phi\circ f_n (x)\,\d x =0\le1.
\]
Therefore, applying (\ref{AA2}) to $f=f_n$ we obtain
\[
P(x,A)=P1_A(x) \le\frac1 n \Phi^{-1} \biggl(\frac1 {\int
_{\mathcal R^d} \mathrm{e}^{-C_\Phi(x,e)}\,\d e} \biggr),
\]
which goes to zero as $n\rightarrow\infty$. Thus $P(x,\cdot)$ is
absolutely continuous w.r.t. the Lebesgue measure, so that the density
function $\varrho (x,y)$ exists, and (\ref{AA1}) follows from
(\ref{AA0}) by taking $f(y)=\Phi^{-1}(\varrho(x,y))$.

Finally, let $\Phi(r)=r^p$ for some $p>1$. For fixed $x$, let
\[
f_n(y)= \frac{ \{n\land\varrho(x, y) \}^{1/(p-1)}}{
(\int_{\mathcal R^{d}} \{
n\land\varrho(x,y) \}^{p/(p-1)}\,\d y )^{1/p}},\qquad n\ge1.
\]
It is easy to see that $\int_{\mathcal R^{d}}f_n^p(y)\,\d y=1$. Then it
follows from (\ref{AA0}) with $\Phi(r)=r^p$ that
\[
\int_{\mathcal R^{d}} \bigl\{n\land\varrho(x,y) \bigr\}^{p/(p-1)}\,\d y \le
\bigl(P f_n(x) \bigr)^{p/(p-1)} \le\frac1 { (\int
_{\mathcal
R^d} \mathrm{e}^{-C_\Phi(x,e)}\,\d e )^{1/(p-1)}}.
\]
Then (\ref{AA}) follows by letting $n\rightarrow\infty$.
\end{pf}

Finally, we consider applications of the shift Harnack inequality to
distribution properties of the underlying transition probability.

\begin{theorem}\label{T4.3} Let $P$ be given above for some Banach space $E$,
and let (\ref{PH}) hold for some $x,e\in E$, finite constant $C_\Phi
(x,e)$ and some strictly increasing and convex continuous function
$\Phi$.
\begin{longlist}[(2)]
\item[(1)] $P(x,\cdot)$ is absolutely continuous w.r.t.
$P(x,\cdot-e)$.

\item[(2)] If $\Phi(r)=r\Psi(r)$ for some strictly increasing
positive continuous function $\Psi$ on $(0,\infty)$. Then the
density $\varrho(x,e;y):=\frac{P(x,\d y)}{P(x,\d y-e)}$
satisfies
\[
\int_E \Phi \bigl(\varrho(x,e;y) \bigr)P(x, \d y-e)\le
\Psi^{-1} \bigl(\mathrm{e}^{C_\Phi
(x,e)} \bigr).
\]
\end{longlist}
\end{theorem}

\begin{pf}
For $P(x,\cdot-e)$-null set $A$, let $f=1_A$. Then (\ref {PH}) implies
that $\Phi(P(x,A))\le0$, hence $P(x,A)=0$ since $\Phi (r)>0$ for $r>0$.
Therefore, $P(x,\cdot)$ is absolutely continuous w.r.t. $P(x,\cdot-e)$.
Next, let $\Phi(r)=r\Psi(r)$. Applying~(\ref{PH}) for $f(y)=
\Psi(n\land\varrho(x,e;y))$ and noting that
\[
Pf(x)=\int_E \bigl\{\Psi \bigl(n\land\varrho(x,e;y)
\bigr) \bigr\}P(x,\d y)\ge \int_E\Phi \bigl(n\land
\varrho(x,e;y) \bigr)P(x,\d y-e),
\]
we obtain
\[
\int_E \Phi \bigl(n\land\varrho(x,e;y) \bigr)P(x, \d y-e)
\le\Psi^{-1} \bigl(\mathrm{e}^{C_\Phi(x,e)} \bigr).
\]
Then the proof is complete by letting $n\rightarrow\infty$.
\end{pf}

\section{Stochastic Hamiltonian system}\label{sec3}

Consider the following degenerate sto\-chastic differential equation on
$\mathcal R^{m+d}=\mathcal R^m\times\mathcal R^d$ $(m\ge0, d\ge1)$:
%
%
\begin{equation}
\label{3.1} \cases{\d X(t)= \bigl\{ AX(t)+BY(t) \bigr\}\,\d t,\vspace*{2pt}
\cr
\d Y(t)= Z \bigl(t,X(t),Y(t) \bigr)\,\d t +\sigma(t)\,\d W(t),}
\end{equation}
where $A$ and $B$ are two matrices of order $m\times m$ and $m\times
d$, respectively, $Z\dvtx  [0,\infty)\times\mathcal R^{m+d}\rightarrow
\mathcal R^d$ is measurable with $Z(t,\cdot)\in C^1(\mathcal R^{m+d})$
for $t\ge0$, $\{\sigma(t)\} _{t\ge0}$ are invertible $d\times
d$-matrices\vspace*{1pt} measurable in $t$ such that the operator norm
$\|\sigma(\cdot)^{-1}\|$ is locally bounded and $W(t)$ is the
$d$-dimensional Brownian motion.

When $m\ge1$ this equation is degenerate, and when $m=0$ we set
$\mathcal R ^m=\{0\}$, so that the first equation disappears and thus,
the equation reduces to a nondegenerate equation on $\mathcal R^d$. To
ensure the existence of the transition density (or heat kernel) of the
associated semigroup $P_t$ w.r.t. the Lebesgue measure on $\mathcal R
^{m+d}$, we make use of the following Kalman rank condition (see
\cite{K}) which implies that the associated diffusion is subelliptic,
{\renewcommand{\theequation}{H}
%
\begin{equation}
\label{eqH} \mbox{There exists $0\le k\le m-1$ such that $\operatorname{Rank}
\bigl[B, AB,\ldots, A^kB \bigr]=m$.}
\end{equation}}\setcounter{equation}{1}%

When $m=0$ this condition is trivial, and for \mbox{$m=1$} it means
that\break  \mbox{$\operatorname{Rank}(B)=1$,} that is, $B\ne0$. For any $m>1$
and $d\ge1$, there exist plenty of examples for matrices $A$ and $B$
such that \textup{(\ref{eqH})} holds; see \cite{K}. Therefore, we allow that
$m$ is much larger than $d$, so that the associated diffusion process
is highly degenerate; see Example~\ref{exa3.1} below.

It is easy to see that if $m=d, \sigma(t)=I_{d\times d}$, $B$ is
symmetric and
\[
Z(x,y)= - \bigl\{\nabla V(x)+ A^*y+F(x,y) (Ax+By) \bigr\}
\]
for some smooth functions $V$ and $F$, then (\ref{3.1}) reduces to the
Hamiltonian system
%
%
\begin{equation}
\label{HS} \qquad\cases{ \d X_t= \nabla H(X_t, \cdot)
(Y_t)\,\d t,\vspace*{2pt}
\cr
\d Y_t =- \bigl\{\nabla
H( \cdot,Y_t) (X_t)+F(X_t,Y_t)
\nabla H(X_t,\cdot) (Y_t) \bigr\}\,\d t +\d W(t)}
\end{equation}
with Hamiltonian function
\[
H(x,y)= V(x)+\langle Ax, y\rangle+\tfrac1 2\langle By,y\rangle;
\]
see, for example, \cite{So}. If, in particular, $A=0,B=I_{d\times d}$
and $F\equiv c$ for some constant $c$, the corresponding Fokker--Planck
equation is known as the ``kinetic Fokker--Planck equation'' in PDE
(see \cite{V}), and the stochastic equation is called ``stochastic
damping Hamiltonian system''; see \cite{Wu}.

Let the solution to (\ref{3.1}) be nonexplosive, and let
\[
P_tf=\mathbf Ef \bigl(X(t),Y(t) \bigr),\qquad t\ge0, f\in
\B_b \bigl(\mathcal R^{m+d} \bigr).
\]
To state our main results, let us fix $T>0$. For nonnegative $\phi\in
C([0,T])$ with $\phi>0$ in $(0,T)$, define
\[
Q_\phi= \int_0^T \phi(t)
\mathrm{e}^{(T-t)A}BB^*\mathrm{e}^{(T-t)A^*}\,\d t.
\]
Then $Q_\phi$ is invertible; cf. \cite{S}. For any $z\in\mathcal
R^{m+d}$ and $r>0$, let $B(z;r)$ be the ball centered at $z$ with
radius $r$.

\begin{theorem}\label{T3.1}
Assume \textup{(\ref{eqH})} and that the solution to (\ref{3.1}) is
nonexplosive such that
%
%
\begin{equation}
\label{UU} \sup_{t\in[0,T]} \mathbf E \Bigl\{\sup
_{B(X(t), Y(t); r)} \bigl|\nabla Z(t,\cdot) \bigr|^2 \Bigr\}<\infty,\qquad
r>0.
\end{equation}
Let $\phi,\psi\in C^1([0,T])$ such that $\phi(0)=\phi(T)=0, \phi
>0$ in
$(0,T)$, and
%
%
\begin{equation}
\label{PS} \psi(T)=1,\qquad \psi(0)=0,\qquad \int_0^T
\psi(t) \mathrm{e}^{(T-t)A}B\,\d t=0.
\end{equation}
Moreover, for $e=(e_1,e_2)\in\mathcal R^{m+d}$, let
\begin{eqnarray*}
h(t) &=& \phi(t) B^*\mathrm{e}^{(T-t)A^*} Q_\phi^{-1}
e_1 +\psi(t) e_2\in\mathcal R^d,
\\
\Theta(t)&=& \biggl(\int_0^t
\mathrm{e}^{(t-s)A} Bh(s)\,\d s, h(t) \biggr)\in\mathcal R ^{m+d},
\qquad t\in[0,T].
\end{eqnarray*}

\begin{longlist}[(2)]
\item[(1)] For any $f\in C_b^1(\mathcal R^{m+d})$, there holds
\begin{eqnarray*}
P_T (\nabla_ef) &=& \mathbf E \biggl\{f \bigl(X(T),Y(T)
\bigr)
\\
&&\hspace*{10pt}{}\times\int_0^T \bigl\langle
\sigma(t)^{-1} \bigl\{h'(t)-\nabla _{\Theta(t)} Z(t,
\cdot) \bigl(X(t),Y(t) \bigr) \bigr\}, \d W(t) \bigr\rangle \biggr\}.
\end{eqnarray*}

\item[(2)] Let $(X(0),Y(0))=(x,y)$ and
\[
R= \exp \biggl[- \int_0^T \bigl\langle
\sigma(t)^{-1}\xi_1(t),\d W(t) \bigr\rangle-
\frac{1} 2 \int_0^T \bigl|
\sigma(t)^{-1}\xi_1(t) \bigr|^2\,\d t \biggr],
\]
where $\xi_1(t)=h'(t)+Z(t,X(t),Y(t))-Z(t,X^1(t),Y^1(t))$ with
\[
X^1(t)=X(t)+\int_0^t \mathrm{e}^{(t-s)A}Bh(s)\,\d s,\qquad Y^1(t)= Y(t) +h(t),\qquad t\ge0.
\]
Then
%
\begin{eqnarray}
\bigl|P_T f(x,y) \bigr|^p &\le& P_T \bigl
\{|f|^p(e+\cdot) \bigr\} (x,y) \bigl(\mathbf ER^{p/(p-1)}
\bigr)^{p-1},\nonumber
\\
\eqntext{p>1, f\in\B _b(E),}
\\
P_T\log f(x,y)&\le&\log P_T \bigl\{f(e+\cdot) \bigr\}(x,y)
+ \mathbf E(R\log R),\nonumber
\\
\eqntext{0<f\in\B_b(E).}
\end{eqnarray}
\end{longlist}
\end{theorem}

\begin{pf} We only prove (1), since (2) follows from Theorem~\ref{T2.1}
with the coupling constructed below for $\varepsilon=1$. Let
$(X^0(t),Y^0(t))=(X(t),Y(t))$ solve (\ref{3.1}) with initial data
$(x,y)$, and for $\varepsilon\in(0,1]$ let $(X^\varepsilon(t),
Y^\varepsilon(t))$ solve the equation
%
%
\begin{equation}
\label{CCC0} %
\cases{ \d X^\varepsilon(t) = \bigl\{A
X^\varepsilon(t) +B Y^\varepsilon (t) \bigr\}\,\d t,
\vspace*{2pt}\cr
\qquad X^\varepsilon(0) =x,
\vspace*{4pt}\cr
\d Y^\varepsilon(t) = \sigma(t)\,
\d W(t) + \bigl\{ Z \bigl(t,X(t),Y(t) \bigr)+\varepsilon h'(t) \bigr
\}\, \d t,
\vspace*{2pt}\cr
\qquad Y^\varepsilon(0)=y.}
\end{equation}
Then it is easy to see that
%
%
\begin{equation}
\label{SOL} \cases{Y^\varepsilon(t) =Y(t) +\varepsilon h(t), \vspace*{2pt}
\cr
\displaystyle X^\varepsilon(t)= X(t) +\varepsilon\int_0^t
\mathrm{e}^{(t-s)A} B h(s)\,\d s.}
\end{equation}
Combining this with $\phi(0)=\phi(T)=0$ and (\ref{PS}), we see that
$h(T)= e_2$ and
\begin{eqnarray*}
&&\int_0^T\mathrm{e}^{(T-t)A} B h(t)\,
\d t
\\
&&\qquad = \int_0^T \phi(t)
\mathrm{e}^{(T-t)A}BB^*\mathrm{e}^{(T-t)A^*}Q_\phi^{-1}
e_1\,\d t +\int_0^T \psi(t)
\mathrm{e}^{(T-t)A} B e_2\,\d t
\\
&&\qquad =e_1.
\end{eqnarray*}
Therefore,
%
%
\begin{equation}
\label{CP} \bigl(X^\varepsilon(T),Y^\varepsilon(T) \bigr)= \bigl(X(T),Y(T)
\bigr)+\varepsilon e,\qquad \varepsilon\in[0,1].
\end{equation}
Next, to see that $ ((X(t),Y(t)), (X^\varepsilon
(t),Y^\varepsilon(t)) )$ is a coupling by change of measure for the
solution to (\ref{3.1}), reformulate (\ref{CCC0}) as
%
%
\begin{equation}
\label{CCC} \qquad\cases{\d X^\varepsilon(t) = \bigl\{A X^\varepsilon(t) +B
Y^\varepsilon(t) \bigr\}\,\d t, &\quad $X^\varepsilon(0) =x$, \vspace*{2pt}
\cr
\d Y^\varepsilon(t) = \sigma(t)\,\d W^\varepsilon(t) + Z
\bigl(t,X^\varepsilon(t),Y^\varepsilon(t) \bigr)\,\d t, &\quad
$Y^\varepsilon(0)=y$,}
\end{equation}
where
%
\begin{eqnarray}
W^\varepsilon(t)&:=& W(t)\nonumber
\\
&&{} +\int_0^t \sigma (s)^{-1}
\bigl\{\varepsilon h'(s)+ Z \bigl(s, X(s),Y(s) \bigr)-Z \bigl(s,
X^\varepsilon(s), Y^\varepsilon (s) \bigr) \bigr\}\,\d s,\nonumber
\\
\eqntext{t\in[0,T].}
\end{eqnarray}
Let
%
%
\begin{equation}
\label{XI} \xi_\varepsilon(s)= \varepsilon h'(s)+ Z \bigl(s,
X(s),Y(s) \bigr)-Z \bigl(s, X^\varepsilon(s), Y^\varepsilon(s) \bigr)
\end{equation}
and
\[
R_\varepsilon=\exp \biggl[- \int_0^T \bigl
\langle\sigma (s)^{-1}\xi_\varepsilon(s),\d W(s) \bigr\rangle -
\frac{1} 2 \int_0^T \bigl|
\sigma(s)^{-1}\xi_\varepsilon(s) \bigr|^2\,\d s \biggr].
\]
By Lemma \ref{LL} below and the Girsanov theorem, $W^\varepsilon(t)$ is
a $d$-dimensional Brownian motion under the probability measure
$\mathbf Q_\varepsilon:=R_\varepsilon\mathbf P$. Therefore, $
((X(t),Y(t)), (X^\varepsilon (t),Y^\varepsilon(t)) )$ is a
coupling by change of measure with changed probability $\mathbf
Q_\varepsilon$. Moreover, combining (\ref{SOL}) with the definition of
$R_\varepsilon$, we see from (\ref{UU}) that
\[
-\frac{\d R_\varepsilon}{\d\varepsilon}\bigg|_{\varepsilon=0}= \int_0^T
\bigl\langle\sigma_s^{-1} \bigl\{ h'(s)-
\nabla_{\Theta(s)} Z (s,\cdot) \bigl(X(s),Y(s) \bigr) \bigr\}, \d W(s) \bigr
\rangle
\]
holds in $L^1(\mathbf P)$. Then the proof is complete by
Theorem~\ref{T2.1}(2).
\end{pf}

\begin{lem}\label{LL} Let the solution to (\ref{3.1}) be nonexplosive
such that (\ref{UU}) holds, and let $\xi_\varepsilon$ be in (\ref{XI}).
Then for any $\varepsilon\in[0,1]$ the process
%
\begin{eqnarray}
R_\varepsilon(t)=\exp \biggl[- \int_0^t \bigl\langle\sigma
(s)^{-1}\xi_\varepsilon(s),\d W(s) \bigr\rangle - \frac{1} 2 \int_0^t
\bigl| \sigma(s)^{-1}\xi_\varepsilon(s) \bigr|^2\,\d s \biggr],\nonumber
\\
\eqntext{t\in[0,T]}
\end{eqnarray}
is a uniformly integrable martingale with $\sup_{t\in[0,T]}\mathbf
E \{ R_\varepsilon(t)\log R_\varepsilon(t) \}<\infty$.
\end{lem}

\begin{pf}
Let $\tau_n=\inf\{t\ge0\dvtx  |X(t)|+|Y(t)|\ge n\}, n\ge1$. Then
$\tau_n\uparrow\infty$ as \mbox{$n\uparrow\infty$}. It suffices to show
that
%
%
\begin{equation}
\label{WWL} \sup_{t\in[0,T],n\ge1}\mathbf E \bigl\{ R_\varepsilon(t
\land \tau _n)\log R_\varepsilon(t\land\tau_n) \bigr
\}< \infty.
\end{equation}
By (\ref{SOL}), there exists $r>0$ such that
%
%
\begin{equation}
\label{QQQ} \bigl(X^\varepsilon(t), Y^\varepsilon (t) \bigr)\in B
\bigl(X(t),Y(t);r \bigr),\qquad t\in [0,T],\varepsilon\in[0,1].
\end{equation}
Let $Q_{\varepsilon,n}= R_\varepsilon(T\land\tau_n)\mathbf P$. By the
Girsanov theorem, $\{ W^\varepsilon(t)\}_{t\in[0,T\land\tau_n]}$ is the
\mbox{$d$-}dimensional Brownian motion under the changed probability
$\mathbf Q_{\varepsilon,n}$. Then, due to (\ref{QQQ}),
\begin{eqnarray*}
&& \sup_{t\in[0,T]}\mathbf E \bigl\{R_\varepsilon(t\land \tau
_n)\log R_\varepsilon(t\land\tau_n) \bigr\}
\\
&&\qquad = \frac1 2 \mathbf E_{\mathbf Q_{\varepsilon,n}} \int_0^{T\land
\tau_n}
\bigl|\sigma(s)^{-1}\xi_\varepsilon(s) \bigr|^2\,\d s
\\
&&\qquad \le C +C \mathbf E_{\mathbf Q_{\varepsilon,n}}\int_0^{T\land\tau
_n}
\sup_{B(X^\varepsilon(t), Y^\varepsilon
(t); r)} \bigl|\nabla Z(t,\cdot) \bigr|^2\,\d t
\end{eqnarray*}
holds for some constant $C>0$ independent of $n$. Since the law of
$(X^\varepsilon(\cdot\land\tau_n), Y^\varepsilon(\cdot\land
\tau_n))$ under $\mathbf Q_{\varepsilon,n}$ coincides with that of
$(X(\cdot\land\tau_n), Y(\cdot\land \tau _n))$ under $\mathbf P$,
combining this with (\ref{UU}), we obtain
\begin{eqnarray*}
&& \sup_{t\in[0,T]}\mathbf E \bigl\{R_\varepsilon(t\land
\tau_n)\log R_\varepsilon(t\land\tau_n) \bigr\}
\\
&&\qquad \le C +C \int_0^{T} \mathbf E\sup
_{B(X(t), Y(t); r)} \bigl|\nabla Z(t,\cdot) \bigr|^2\,\d t<\infty.
\end{eqnarray*}
Therefore, (\ref{WWL}) holds.
\end{pf}

\begin{rem}\label{rem3.1}
(a) As shown in \cite{GW}, Lemma~2.4, condition (\ref{UU}) is implied
by the Lyapunov condition {(A)} therein, for which some concrete
examples have been presented in \cite{GW}. Moreover, as shown in
\cite{GW}, Section~3 (see also Theorem~4.1 in~\cite{WZ}) that under
reasonable grown conditions of $\nabla Z(t,\cdot)$ one obtains from
Theorem \ref{T3.1}(1)
%
\begin{eqnarray}
P_t |\nabla f|\le\delta \bigl\{P_t(f\log f)-(P_t f)\log P_t f
\bigr\} +\frac{W(t,\cdot)}\delta P_t f,\nonumber
\\
\eqntext{ t>0, f\in\B _b^+
\bigl(\mathcal R^{m+d} \bigr), \delta>\delta_0}
\end{eqnarray}
for some constant $\delta_0\ge0$ and some positive functions $W(t,\cdot
)$. According to Theorem \ref{MA}, this inequality implies the shift
Harnack inequality.

(b) For any $T_2>T_1$. Applying Theorem \ref{T3.1} to $ (X(T_1+t),
Y(T_1+t))$ in place of $(X(t),Y(t))$, we see that the assertions in
Theorem \ref{T3.1} hold for
\[
P_{T_1,T_2} f(x,y):=\mathbf E \bigl(f \bigl(X(T_2),
Y(T_2) \bigr) \mid \bigl(X(T_1),Y(T_1)
\bigr)=(x,y) \bigr)
\]
in place of $P_Tf$ with $T$ and $0$ replaced by $T_2$ and $T_1$,
respectively.
\end{rem}

To derive explicit inequalities from Theorem \ref{T3.1}, we consider
below a special case where $\|\nabla Z(t,\cdot)\|_\infty$ is bounded
and $A^l=0$ for some natural number $l\ge1$.

\begin{cor}\label{C3.2}
Assume \textup{(\ref{eqH})}. If $\|\nabla Z(t,\cdot)\| _\infty$ and
$\|\sigma(t)^{-1}\|$ are bounded in $t\ge0$, and $A^l=0$ for some
$l\ge1$. Then there exists a constant $C>0$ such that for any positive
$f\in\B_b(\mathcal R^{m+d}), T>0$ and $e=(e_1,e_2)\in\mathcal R^{m+d}$:
\begin{longlist}[(2)]
\item[(1)] $(P_Tf)^p \le P_T \{f ^p(e+\cdot)\} \exp [\frac
{Cp}{p-1} (\frac{|e_2|^2}{1\land T}
+\frac{|e_1|^2}{(1\land T)^{4k+3}} ) ]$, $p>1$;

\item[(2)] $ P_T\log f \le\log P_T \{f(e+\cdot)\} + C (\frac
{|e_2|^2}{1\land T} +\frac{|e_1|^2}{(1\land T)^{4k+3}} )$;

\item[(3)] for $f\in C_b^1(\mathcal R^{m+d})$, $| P_T \nabla_e f|^2
\le C |P_T f^2|  (\frac{|e_2|^2}{1\land T}
+\frac{|e_1|^2}{(1\land T)^{4k+3}} )$;

\item[(4)] for strictly positive $f\in C_b^1(\mathcal R^{m+d})$,
\begin{eqnarray*}
\bigl|P_T \nabla_e f \bigr|(x,y)&\le&\delta \bigl \{P_T(f\log
f)-(P_T f)\log P_T f \bigr\}
\\
&&{} + \frac{C}\delta
\biggl(\frac{|e_2|^2}{1\land T} +\frac {|e_1|^2}{(1\land T)^{4k+3}}
\biggr) P_T f,\qquad \delta>0.
\end{eqnarray*}
\end{longlist}
\end{cor}

\begin{pf}
According to Remark~\ref{rem3.1}(b), $P_T= P_{T-1}P_{T-1,T}$ and the
Jensen inequality, we only need to prove\vspace*{1pt} for $T\in(0,1]$. Let $\phi
(t)= \frac{t(T-t)}{T^2}$. Then $\phi(0)=\phi(T)=0$ and due to
\cite{WZ}, Theorem~4.2(1), the rank condition \textup{(\ref{eqH})} implies that
%
%
\begin{equation}
\label{Q} \bigl\|Q_\phi^{-1} \bigr\|\le c T^{-(2k+1)}
\end{equation}
for some constant $c>0$ independent of $T\in(0,1]$. To fix the other
reference function $\psi$ in Theorem \ref{T3.1}, let $\{c_i\}_{1\le
i\le l+1}\in\mathcal R$ be such that
\[
\cases{ \displaystyle 1+\sum_{i=1}^{l+1}c_i=0,
\vspace*{2pt}
\cr
\displaystyle 1+\sum_{i=1}^{l+1}
\frac{j+1}{j+1+i} c_i=0, &\quad $0\le j\le l-1$.}
\]
Take
\[
\psi(t)= 1+\sum_{i=1}^{l+1} c_i
\frac{(T-t)^i}{T^i},\qquad t\in[0,T].
\]
Then $\psi(0)=0,\psi(T)=1$ and $\int_0^T (T-t)^j\psi(t)\,\d t=0$ for
$0\le j\le l-1$. Since $A^l=0$, we conclude that $\int_0^T\psi
(t)\mathrm{e} ^{(T-t)A}\,\d t=0$. Therefore, (\ref{PS}) holds. It is
easy to see that
\[
\bigl|\psi(t)\bigr| \le c,\qquad \bigl|\psi'(t)\bigr|\le c T^{-1},\qquad t
\in[0,T]
\]
holds for some constant $c>0$. Combining this with (\ref{Q}), (\ref
{SOL}) and the boundedness of $\|\nabla Z\|_\infty$ and $\|\sigma
^{-1}\|$, we obtain
%
%
\begin{eqnarray}
\label{**} \bigl|\xi_1(t) \bigr|+ \bigl|h'(t) \bigr|&\le& c
\bigl(T^{-2(k+1)}|e_1|+T^{-1}|e_2| \bigr),
\nonumber
\\[-8pt]
\\[-8pt]
\bigl|\Theta(t) \bigr|&\le& c \bigl(T^{-(2k+1)}|e_1|+ |e_2|
\bigr)
\nonumber
\end{eqnarray}
for some constant $c>0$. From this and Theorem \ref{T3.1}, we derive
the desired assertions.
\end{pf}

\begin{cor}
In the situation of Corollary~\ref{C3.2}. Let $\|\cdot\| _{p\rightarrow
q}$ be the operator norm from $L^p$ to $L^q$ w.r.t. the Lebesgue
measure on $\mathcal R^{m+d}$. Then there exists a constant $C>0$ such
that
%
%
\begin{eqnarray}
\label{CC} \|P_T\|_{p\rightarrow\infty}\le C^{1/p} \biggl(\frac
p{p-1} \biggr)^{(m+d)/(2p)} (1\land T)^{-(d+(4k+3)m)/(2p)},
\nonumber\\[-12pt]\\[-12pt]
\eqntext{p>1, T>0.}
\end{eqnarray}
Consequently, the transition density $p_T((x,y),(x',y'))$ of $P_T$
w.r.t. the Lebesgue measure on $\mathcal R^{m+d}$ satisfies
%
%
\begin{eqnarray}
\label{HHH}
\qquad && \int_{\mathcal R^{m+d}} p_T \bigl((x,y),
\bigl(x',y' \bigr) \bigr)^{p/(p-1) }\,\d
x'\,\d y'\nonumber
\\
&&\qquad \le C^{1/(p-1)} \biggl(\frac p{p-1} \biggr)^{(m+d)/(2(p-1))} (1\land
T)^{-(d+(4k+3)m)/(2(p-1))},
\\
\eqntext{T>0, (x,y)\in\mathcal R^{m+d},p>1.}
\end{eqnarray}
\end{cor}

\begin{pf}
By Corollary~\ref{C3.2}(1), (\ref{CC}) follows from (\ref {AA0}) for
$P_T=P$, $\Phi(r)=r^p$ and
\[
C_\Phi \bigl((x,y),(e_1,e_2) \bigr)=
\frac{Cp}{p-1} \biggl(\frac{|e_2|^2}{1\land
T} +\frac
{|e_1|^2}{(1\land T)^{4k+3}} \biggr).
\]
Moreover, (\ref{HHH}) follows from (\ref{AA}).
\end{pf}

\begin{exa}\label{exa3.1}
A simple example for Theorem \ref{C3.2} to
hold is that $\sigma(t)=\sigma$ and $Z(t,\cdot)=Z$ are independent of
$t$ with $\|\nabla Z\|_\infty<\infty$, $A=0$ and\break
$\operatorname{Rank}(B)=m$. In this case we have $d\ge m$; that is, the
dimension of the generate part is controlled by that of the
nondegenerate part. In general, our results allow $m$ to be much
larger\vadjust{\goodbreak}
than $d$. For instance, let $m=ld$ for some $l\ge2$ and
\[
A=\pmatrix{0 &I_{d\times d} &0 &\cdots&0 &0
\cr
0 &0 &I_{d\times d} &
\cdots&0 &0
\cr
\cdots&\cdots&\cdots&\cdots&\cdots&\cdots
\cr
0 &0 &0 &\cdots&0
&I_{d\times d}
\cr
0 &0 &0 &\cdots&0 &0}_{(ld)\times(ld)},\qquad B= \pmatrix{0
\cr
\cdot
\cr
\cdot
\cr
\cdot
\cr
0
\cr
I_{d\times d}}_{(ld)\times d}.
\]
Then $A^l=0$ and \textup{(\ref{eqH})} holds for $k=m-1$. Therefore, assertions
in Corollary~\ref{C3.2} hold for $k=l-1$.
\end{exa}

\section{Functional stochastic differential equations}\label{sec4}

The purpose of this section is to establish Driver's integration by
parts formula and shift Harnack inequality for delayed stochastic
differential equations. In this case the associated segment processes
are functional-valued, and thus, infinite-dimensional. As continuation
to Section~\ref{sec3}, it is natural for us to study the generalized
stochastic Hamiltonian system with delay as in \cite{BWY}, where the
Bismut formula and the Harnack inequalities are derived using coupling.
However, for this model it seems very hard to construct the required
new-type couplings. So, we only consider here the nondegenerate
setting.

Let $\tau>0$ be a fixed number, and let $\mathcal C=C([-\tau,0];\mathcal R^d)$ be equipped with uniform norm $\|\cdot\|_\infty$.
For simplicity, we will use $\nabla$ to denote the gradient operator
both on $\mathcal R^d$ and $\mathcal C$. For instance, for a
differentiable function $F$ on $\mathcal C$ and $\xi \in\mathcal C$,
$\nabla F(\xi)$ is a linear operator from $\mathcal C$ to $\mathcal R$
with
\[
\mathcal C\ni\eta\mapsto\nabla_\eta F(\xi)=\lim_{\varepsilon
\rightarrow0}
\frac{F(\xi+\varepsilon\eta
)-F(\xi)}\varepsilon.
\]
Moreover, let $\|\cdot\|$ be the operator norm for linear operators.
Finally, for a function $h\in C([-\tau,\infty);\mathcal R^d)$ and
$t\ge0$, let $h_t\in\mathcal C$ be such that
$h_t(\theta)=h(t+\theta)$,
$\theta\in [-\tau,0]$.

Consider the following stochastic differential equations on $\mathcal
R^d$:
%
%
\begin{equation}
\label{4.1'} \d X(t)=b(t,X_t)\,\d t +\sigma(t)\,\d
W(t),\qquad t\ge0,
\end{equation}
where $W(t)$ is the Brownian motion on $\mathcal R^d$, $b\dvtx
[0,\infty )\times\mathcal C\rightarrow \mathcal R^d$ is measurable such
that $\|\nabla b(t,\cdot)\|_\infty$ is locally bounded in $t$, and
$\sigma\dvtx  [0,\infty)\rightarrow\mathcal R^d\otimes\mathcal R^d$ is
measurable with $\|\sigma(t)^{-1}\|$ locally bounded. We remark that
the local boundedness assumption of $\|\nabla b(t,\cdot)\|_\infty$ is
made only for simplicity and can be weakened by some growth conditions
as in \cite{BWY}.

Now, for any $\xi\in\mathcal C$, let $X^\xi(t)$ be the solution to
(\ref{4.1'}) for $X_0=\xi$, and let $X_t^\xi$ be the associated segment
process. Let
\[
P_tF(\xi)=\mathbf EF \bigl(X_t^\xi \bigr),
\qquad t\ge0, \xi\in\mathcal C, F\in \B_b(\mathcal C).
\]
We aim to establish the integration by parts formula and shift Harnack
inequality for $P_T$. It turns out that we are only\vadjust{\goodbreak} able to make
derivatives or shifts along directions in the Cameron--Martin space
\[
\mathcal H:= \biggl\{h\in\mathcal C\dvtx \|h\|_\mathcal H^2:=
\int_{-\tau}^0 \bigl|h'(t)\bigr|^2 \,\d t <\infty \biggr\}.
\]

\begin{theorem}\label{T3.1'} Let $T>\tau$ and $\eta\in\mathcal H$ be
fixed. For any $\phi\in\B_b([0,T-\tau])$ such that $ \int_0^{T-\tau}
\phi(t)\,\d t =1$, let
\[
\Gamma(t)= \cases{\phi(t)\eta(-\tau), &\quad if $t\in[0,T-\tau]$, \vspace*{2pt}
\cr
\eta'(t-T), &\quad if $t\in(T-\tau,T]$.}
\]
Let $\|\sigma(t)^{-1}\|\le K(T), \|\nabla b(t,\cdot)\|_\infty\le
\kappa(T)$ for $t\in[0,T]$.
\begin{longlist}[(2)]
\item[(1)] For any $F\in C_b^1(\mathcal C)$,
\[
P_T(\nabla_\eta F) = \mathbf E \biggl(F(X_T)
\int_0^T \bigl\langle \sigma(t)^{-1}
\bigl(\Gamma (t)-\nabla_{\Theta_t}b(t,\cdot) (X_t) \bigr), \d
W(t) \bigr\rangle \biggr)
\]
holds for
\[
\Theta(t)= \int_0^{t\lor0} \Gamma(s)\,\d s,\qquad t
\in[-\tau, T].
\]
Consequently, for any $\delta>0$ and positive $F\in C_b^1(\mathcal C)$,
\begin{eqnarray*}
\bigl|P_T(\nabla_\eta F) \bigr|&\le&\delta \bigl\{P_T(F
\log F)-(P_TF)\log P_T F \bigr\}
\\
&&{}  +\frac{2K(T)^2(1+\kappa(T)^2T^2)} \delta
\biggl(\|\eta\|_\mathcal H^2 +\frac{|\eta(-\tau)|^2}{T-\tau}
\biggr)P_TF.
\end{eqnarray*}

\item[(2)] For any nonnegative $F\in\B_b(\mathcal C)$,
\begin{eqnarray*}
(P_TF)^p &\le& \bigl(P_T \bigl\{F(\eta+ \cdot)
\bigr\}^p \bigr)
\\
&&{}\times  \exp \biggl[\frac{2pK(T)^2(1+\kappa(T)^2T^2)} {p-1} \biggl(\| \eta
\|_\mathcal H^2 +\frac{|\eta(-\tau)|^2}{T-\tau} \biggr) \biggr].
\end{eqnarray*}

\item[(3)] For any positive $F\in\B_b(\mathcal C)$,
\[
P_T\log F \le\log P_T \bigl\{F(\eta+\cdot) \bigr\}+ 2
K(T)^2 \bigl(1+\kappa(T)^2T^2 \bigr) \biggl(
\| \eta\|_\mathcal H^2 +\frac{|\eta(-\tau)|^2}{T-\tau} \biggr).
\]
\end{longlist}
\end{theorem}

\begin{pf}
For fixed $\xi\in\mathcal C$, let $X(t)$ solve (\ref{4.1'}) for
$X_0=\xi$. For any $\varepsilon\in[0,1]$, let $X^\varepsilon(t)$ solve
the equation
\[
\d X^\varepsilon(t)= \bigl\{b(t,X_t)+\varepsilon\Gamma(t) \bigr
\}\, \d t +\sigma(t)\,\d W(t),\qquad t\ge0, X_0^\varepsilon=\xi.
\]
Then it is easy to see that
%
%
\begin{equation}
\label{TH0} X_t^\varepsilon= X_t+\varepsilon\Theta
_t,\qquad t\in[0,T].
\end{equation}
In particular, $X_T^\varepsilon=X_T+\varepsilon\eta$. Next, let
\begin{eqnarray*}
R_\varepsilon&=& \exp \biggl[ -\int_0^T
\bigl\langle\sigma (t)^{-1} \bigl\{\varepsilon\Gamma(t)+b(t,X_t)-b
\bigl(t,X_t^\varepsilon \bigr) \bigr\}, \d W(t) \bigr\rangle
\\
&&\hspace*{26pt}{} -\frac{1}{2} \int_0^T \bigl|
\sigma(t)^{-1} \bigl\{\varepsilon \Gamma(t)+b(t,X_t)-b
\bigl(t,X_t^\varepsilon \bigr) \bigr\} \bigr|^2\,\d t
\biggr].
\end{eqnarray*}
By the Girsanov theorem, under the changed probability $\mathbf
Q_\varepsilon:=R_\varepsilon\mathbf P$, the process
\[
W^\varepsilon(t):= W(t)+\int_0^t
\sigma(s)^{-1} \bigl(\Gamma (s)+b(s,X_s)-b
\bigl(s,X_s^\varepsilon \bigr) \bigr)\,\d s,\qquad t\in[0,T]
\]
is a $d$-dimensional Brownian motion. So, $(X_t,X^\varepsilon_t)$ is a
coupling by change of measure with changed probability $\mathbf
Q_\varepsilon $. Then the desired integration by parts formula follows
from Theorem \ref{T2.1} since $R_0=1$ and due to (\ref{TH0}),
\[
\frac{\d}{\d\varepsilon}R^\varepsilon\bigg|_{\varepsilon=0} = -\int
_0^T \bigl\langle\sigma(t)^{-1}
\bigl(\Gamma (t)-\nabla_{\Theta_t}b(t,\cdot) (X_t) \bigr), \d
W(t) \bigr\rangle
\]
holds in $L^1(\mathbf P)$. Taking $\phi(t)= \frac1{T-\tau}$, we have
\begin{eqnarray*}
\int_0^T \bigl|\Gamma(t) \bigr|^2\,\d t &
\le&\|\eta\|_\mathcal H^2+ \frac{|\eta(-\tau)|^2}{T-\tau},
\\
\bigl\|\nabla_{\Theta_t}b(t,\cdot) \bigr\|_\infty^2&\le&
\kappa(T)^2 \biggl(\int_0^T \bigl|
\Gamma(t) \bigr|\,\d t \biggr)^2\le\kappa(T)^2T \int
_0^T \bigl|\Gamma(t) \bigr|^2\,\d t.
\end{eqnarray*}
Then
%
%
\begin{equation}
\label{LLO} \qquad\quad\int_0^T \bigl|\Gamma(t)-
\nabla_{\theta
_t}b(t,\cdot) (X_t) \bigr|^2\,\d t \le2
\bigl(1+T^2\kappa(T)^2 \bigr) \biggl(\|\eta
\|_\mathcal H^2 +\frac{|\eta(-\tau)|^2}{T-\tau} \biggr).
\end{equation}
So,
\begin{eqnarray*}
&&\log\mathbf E\exp \biggl[ \frac1 \delta\int_0^T
\bigl\langle\sigma(t)^{-1} \bigl\{\Gamma(t)-\nabla_{\Theta
_t}b(t,
\cdot) (X_t) \bigr\}, \d W(t) \bigr\rangle \biggr]
\\
&&\qquad\le\frac1 2 \log\mathbf E\exp \biggl[ \frac{2 K(T)^2}{\delta
^2} \int
_0^T \bigl|\Gamma (t)-\nabla_{\Theta_t}b(t,
\cdot) (X_t) \bigr|^2\,\d t \biggr]
\\
&&\qquad\le\frac{2K(T)^2(1+T^2\kappa(T)^2)}{\delta^2} \biggl(\|\eta\| _\mathcal H^2 +
\frac{|\eta(-\tau)|^2}{T-\tau} \biggr).
\end{eqnarray*}
Then the second result in (1) follows from the Young inequality
\begin{eqnarray*}
\bigl|P_T(\nabla_\eta F) \bigr|&\le& \delta \bigl
\{P_T(F\log F)-(P_TF)\log P_T F \bigr\}
\\
&&{}+\delta\log\mathbf E\exp \biggl[ \frac1 \delta\int_0^T
\bigl\langle\sigma(t)^{-1} \bigl\{\Gamma (t)-\nabla_{\Theta_t}b(t,
\cdot) (X_t) \bigr\}, \d W(t) \bigr\rangle \biggr].
\end{eqnarray*}
Finally, (2) and (3) can be easily derived by applying Theorem
\ref{T2.1} for the above constructed coupling with $\varepsilon=1$, and
using (\ref {TH0}) and (\ref{LLO}).
\end{pf}

From Theorem \ref{T3.1'} we may easily derive regularization estimates
on $P_T(\xi,\cdot)$, the distribution of $X^\xi_T$. For instance,
Theorem \ref{T3.1'}(1) implies estimates on the derivative of
$P_T(\xi,A+\cdot)$ along $\eta\in\mathcal H$ for $\xi\in\mathcal C$ and
measurable $A\subset\mathcal C $; and due to Theorems \ref{T4.3},
\ref{T3.1'}(2) and \ref{T3.1'}(3) imply some integral estimates on the
density $p_T(\xi,\eta;\gamma):= \frac{P_T(\xi,\d\gamma)}{P_T(\xi,\d
\gamma-\eta)}$ for $\eta \in\mathcal H$. Moreover, since $\mathcal H$
is dense in $\mathcal C$, the shift Harnack inequality in Theorem
\ref{T3.1'}(2) implies that $P_T(\xi,\cdot)$ has full support on
$\mathcal C$ for any $T>\tau$ and $\xi\in\mathcal C$.

\section{Semi-linear stochastic partial differential equations}\label{sec5}

The purpose of this section is to establish Driver's integration by
parts formula and shift Harnack inequality for semi-linear stochastic
partial differential equations. We note that the Bismut formula has
been established in \cite{BWY2} for a class of delayed SPDEs, but for
technical reasons we only consider here the case without delay.

Let $(H,\langle\cdot,\cdot\rangle,|\cdot|)$ be a real separable Hilbert
space, and $(W(t))_{t\geq0}$ a cylindrical Wiener process on $H$ with
respect to a complete probability space $(\Omega, \mathcal{F},
\mathbf{P})$ with the natural filtration $\{\mathcal{F}_t\}_{t\geq0}$.
Let $\mathcal{L}(H)$ and $\mathcal {L}_{\mathrm{HS}}(H)$ be the spaces of all
linear bounded operators and Hilbert--Schmidt operators on $H$,
respectively. Denote by $\|\cdot\|$ and $\|\cdot\|_{\mathrm{HS}}$ the operator
norm and the Hilbert--Schmidt norm, respectively.

Consider the following semi-linear SPDE:
%
%
%
\begin{equation}
\label{eq1} \cases{ \d X(t)= \bigl\{AX(t)+b \bigl(t,X(t) \bigr) \bigr\}\,\d t+
\sigma(t)\,\d W(t), \vspace*{2pt}
\cr
X(0)=x\in H,}
\end{equation}
where
\begin{longlist}[(A1)]
\item[(A1)] $(A,\mathcal{D}(A))$ is a linear operator on $H$
generating a contractive, strongly continuous semigroup
$(\mathrm{e}^{tA})_{t\geq0}$ such that $\int_0^1
\|\mathrm{e}^{sA}\| _{\mathrm{HS}}^2\,\d s<\infty$.

\item[(A2)] $b\dvtx [0,\infty)\times H \rightarrow H$ is measurable,
and Fr\'{e}chet differentiable in the second variable such that
$\|\nabla b(t,\cdot)\|_\infty:= \sup_{x\in H}\|\nabla b(t,\cdot
)(x)\|$ is locally bounded in $t\ge0$.

\item[(A3)] $ \sigma\dvtx  [0,\infty)\rightarrow\mathcal{L}(H)$ is
measurable and locally bounded, and $\sigma(t)$ is invertible such
that $\| \sigma(t)^{-1}\|$ is locally bounded in $t\ge0$.
\end{longlist}
Then the equation (\ref{eq1}) has a unique a mild solution (see
\cite{DZ}), which is an adapt process $(X(t))_{t\ge0}$ on $H$ such that
\[
X(t)= \mathrm{e}^{t A}x + \int_0^t
\mathrm{e}^{(t-s)A} b \bigl(s,X(s) \bigr)\,\d s +\int_0^t
\mathrm{e} ^{(t-s)A}\sigma(s)\,\d W(s),\qquad t\ge0.
\]
Let
\[
P_tf \bigl(X(0) \bigr) =\mathbf{E}f \bigl(X(t) \bigr),\qquad t
\geq0, X(0)\in H, f\in\B_b(H).
\]
Finally, for any $e\in H$, let
\[
e(t)= \int_0^t\mathrm{e}^{sA} e\,
\d s,\qquad t\ge0.
\]

\begin{theorem}\label{T5.1}
Let $T>0$ and $e\in\mathcal D(A)$ be fixed. Let $\|\sigma (t)^{-1}\|\le
K(T)$, $\|\nabla b(t,\cdot)\|_\infty\le\kappa(T)$ for $t\in[0,T]$.
\begin{longlist}[(2)]
\item[(1)] For any $f\in C_b^1(H)$,
\[
P_T(\nabla_{e(T)} f) = \mathbf E \biggl(f \bigl(X(T) \bigr)
\int_0^T \bigl\langle \sigma(t)^{-1}
\bigl(e-\nabla _{e(t)}b(t,\cdot) \bigl(X(t) \bigr) \bigr), \d W(t) \bigr
\rangle \biggr).
\]
Consequently, for any $\delta>0$ and positive $f\in C_b^1(H)$,
\begin{eqnarray*}
\bigl|P_T(\nabla_{e(T)} f) \bigr|&\le&\delta \bigl\{P_T(f
\log f)-(P_Tf)\log P_T f \bigr\}
\\
&&{} +\frac{K(T)^2|e|^2}{\delta}
\biggl( T +T^2\kappa(T)+\frac
{T^3\kappa(T)^2}3 \biggr)P_Tf.
\end{eqnarray*}

\item[(2)] For any nonnegative $F\in\B_b(H)$,
\[
(P_TF)^p \le \bigl(P_T \bigl\{F \bigl(e(T)+
\cdot \bigr) \bigr\}^p \bigr) \exp \biggl[ \frac{pK(T)^2|e|^2}{p-1} \biggl(
T +T^2\kappa (T)+\frac{T^3\kappa
(T)^2}3 \biggr) \biggr].
\]

\item[(3)] For any positive $F\in\B_b(H)$,
\begin{eqnarray*}
P_T\log F &\le& \log P_T \bigl\{F \bigl(e(T)+\cdot \bigr)
\bigr\}
\\
&&{} + K(T)^2|e|^2 \biggl( T +T^2\kappa (T)+
\frac{T^3\kappa(T)^2}3 \biggr).
\end{eqnarray*}
\end{longlist}
\end{theorem}

\begin{pf}
For fixed $x\in H$, let $X(t)$ solve (\ref{4.1'}) for $X(0)=x$. For any
$\varepsilon\in[0,1]$, let $X^\varepsilon(t)$ solve the equation
%
\begin{eqnarray}
\d X^\varepsilon(t)= \bigl\{AX^\varepsilon(t)+b \bigl(t,X(t) \bigr)+
\varepsilon e \bigr\}\,\d t +\sigma(t)\,\d W(t),\nonumber
\\
\eqntext{t\ge0, X^\varepsilon(0)=x.}
\end{eqnarray}
Then it is easy to see that
%
%
\begin{equation}
\label{TH} X^\varepsilon(t) = X(t)+\varepsilon e(t),\qquad t\in[0,T].
\end{equation}
In particular, $X^\varepsilon(T) =X(T)+\varepsilon e(T)$. Next, let
\begin{eqnarray*}
R_\varepsilon &=& \exp \biggl[ -\int_0^T
\bigl\langle\sigma (t)^{-1} \bigl\{\varepsilon e +b \bigl(t,X(t)
\bigr)-b \bigl(t,X^\varepsilon(t) \bigr) \bigr\}, \d W(t) \bigr\rangle
\\
&&\hspace*{26pt}{}-\frac1
2 \int_0^T \bigl|\sigma(t)^{-1} \bigl\{
\varepsilon e +b \bigl(t,X(t) \bigr)-b \bigl(t,X^\varepsilon (t) \bigr) \bigr\}
\bigr|^2\,\d t \biggr].
\end{eqnarray*}
By the Girsanov theorem, under the weighted probability $\mathbf
Q_\varepsilon:=R_\varepsilon\mathbf P$, the process
\[
W^\varepsilon(t):= W(t)+\int_0^t
\sigma(s)^{-1} \bigl(\varepsilon e+b(s,X_s)-b
\bigl(s,X_s^\varepsilon \bigr) \bigr)\,\d s,\qquad t\in[0,T]
\]
is a $d$-dimensional Brownian motion. So, $(X(t),X^\varepsilon(t))$ is
a coupling by change of measure with changed probability $\mathbf
Q_\varepsilon$. Then the desired integration by parts formula follows
from Theorem \ref {T2.1} since $R_0=1$ and due to (\ref{TH}),
\[
\frac{\d}{\d\varepsilon}R^\varepsilon\bigg|_{\varepsilon=0} = -\int_0^T
\bigl\langle\sigma(t)^{-1} \bigl(e-\nabla _{e(t)}b(t,\cdot)
\bigl(X(t) \bigr) \bigr), \d W(t) \bigr\rangle
\]
holds in $L^1(\mathbf P)$. This formula implies the second inequality
in (1) due to the given upper bounds on $\|\sigma(t)^{-1}\|$ and
$\|\nabla b(t,\cdot )\|$ and the fact that
\begin{eqnarray*}
&& \bigl|P_T(\nabla_\eta F) \bigr|- \delta \bigl\{P_T(F
\log F)-(P_TF)\log P_T F \bigr\}
\\
&&\qquad \le\delta\log\mathbf E\exp \biggl[ \frac1 \delta\int
_0^T \bigl\langle\sigma(t)^{-1}
\bigl(e-\nabla _{e(t)}b(t,\cdot) \bigl(X(t) \bigr) \bigr), \d W(t) \bigr
\rangle \biggr]P_TF
\\
&&\qquad \le\frac\delta2 \log\mathbf E\exp \biggl[\frac2 {\delta^2}
\int_0^T \bigl|\sigma (t)^{-1} \bigl(e-
\nabla_{e(t)}b(t,\cdot) \bigl(X(t) \bigr) \bigr) \bigr|^2\,\d t
\biggr]P_TF.
\end{eqnarray*}
Finally, since $|e(t)|\le t|e|$, (2) and (3) can be easily derived by
applying Theorem~\ref{T2.1} for the above constructed coupling with
$\varepsilon=1$.
\end{pf}

\section*{Acknowledgment}
The author would like to thank the referees for helpful comments and corrections.



\printaddresses

\end{document}